\documentclass[11pt,leqno]{article}
\usepackage{amsfonts}
\usepackage{amsmath}
\def\CC {{\mathbb C}}

\def\NN {{\mathbb N}}
\def\ZZ {{\mathbb Z}}

\def\CP {{\mathbb C}{\mathbb P}^{1}}
\def\be {\begin{eqnarray}}
\def\ben {\begin{eqnarray*}}
\def\ee {\end{eqnarray}}
\def\een {\end{eqnarray*}}

\def\AAA{\kern-0.3em}
\def\AA{\kern-0.18em}
\def\AC{\kern-0.14em}
\def\AB{\kern-0.22em}

\newtheorem{th1}{Theorem}

\title{ Non-integrability of the second Painlev\'e
               equation as a Hamiltonian system}
\author{Tsvetana Stoyanova, Ognyan Christov}

\date{}

\begin{document}
 \maketitle

{\bf Abstract}

     We prove non-existence of an additional rational integral for
     the second Painlev\'e equation
     $({\rm P}_{{\rm II}})$  considered as a Hamiltonian system using
     Morales - Ramis theory.

{\bf Key words: Painlev\'e second equation, Differential Galois Theory, Hamiltonian system}

{\bf 2000 Mathematics Subject Classification: 70H07, 34M55, 37J30}

\section{Introduction }

   The six Painlev\'e equations ${\rm P}_{{\rm I}} - {\rm P}_{{\rm VI}}$
   are non - linear second
   order differential equations with the property that their
   solutions have no movable critical points $[^1]$.
   Except ${\rm P}_{{\rm I}}$ the equations ${\rm P}_{{\rm II}} - {\rm P}_{{\rm VI}}$
   depend on
   some parameters. These six equations are not reducible to
   linear equations, nor integrable in terms of previously known
   functions. Hence they define new transcendental functions for their
   general solutions called Painlev\'e transcendents.

   Nowadays  Painlev\'e transcendents have important applications in the
   mathematical physics. We only mention that they arise as
   reductions of soliton equations such as KdV, nonlinear
   Schr\"odinger,  Kadomtsev-Petviashvili equation, etc.

   Three essential  properties of the Painlev\'e equations are the
   following $[^{1,2}]$:
    \begin{enumerate}

            \vspace{-1ex}
            \item

     They admit a Hamiltonian formulation;

            \vspace{-1ex}
            \item
     They can be expressed as isomonodromic deformation of some
     linear differential equations with rational coefficients;

             \vspace{-1ex}
             \item
     Except for ${\rm P}_{{\rm I}}$ all the other five equations admit
     one-parameter families of solutions by means of special
     functions and also for some  special values of the
     parameters they have particular rational solutions.

     \end{enumerate}

      As the Painlev\'e equations are Hamiltonian systems, then
      their integrability should be considered in the context of
      Hamiltonian systems.
       In $[^2]$ Morales-Ruiz asked the question to prove rigorously
       non-integrability of the Hamiltonian systems
       equivalent to Painlev\'e equations which have rational
       particular solutions. In this note we consider  ${\rm P}_{{\rm II}}$.

       The second Painlev\'e equation $({\rm P}_{{\rm II}})$
           \be
           \label{pe}
               \frac{d^2w}{dz^2}=2 w^3 + z w + \alpha ,
          \ee
       where $\alpha$ is a parameter, can be written in the
       Hamiltonian form
       $
                \frac{dq}{dz}  =  \frac{\partial H}{\partial p}, \quad
                \frac{dp}{dz}  =
               - \frac{\partial H}{\partial q}
               $
       with the Hamiltonian
          $    H(q, p, \alpha) =
                \frac{1}{2}p^2 - (q^2 + \frac{1}{2}z) p - (\alpha + \frac{1}{2})\,q.$
       This system  is non-autonomous. Next, we
        extent it to a two degrees of freedom Hamiltonian system.
       We introduce the conjugate variable $E$ to $z$ in the
       standard way $[^4]$, namely putting
       $H = \frac{1}{2}p^2 - (q^2 + \frac{1}{2}z)p - (\alpha + \frac{1}{2})\,q - E = 0.$
       Denote $F:= -E.$ Then we obtain a Hamiltonian
           \be
              \label{ham}
               \widehat{H}(q, p, z, F)=
                 \frac{1}{2}p^2 + F - (q^2 + \frac{1}{2}z)p - (\alpha + \frac{1}{2})\,q
           \ee
       defined on $M := \big\{(q, p, z, F) \in \CC^4\big\}$ with the
       canonical symplectic structure $dp\wedge dq + dF\wedge dz.$
       The corresponding autonomous Hamiltonian system reads
          \be
            \label{eqh}
               \frac{dq}{ds} = \frac{\partial \widehat{H}}{\partial p}
                              = p - q^2-\frac{1}{2}z , \qquad \qquad
                      \frac{dz}{ds} = \frac{\partial \widehat{H}}{\partial F}
                              = 1            \\
               \frac{dp}{ds} = - \frac{\partial \widehat{H}}{\partial q}
                              =  2qp + \alpha + \frac{1}{2},  \qquad
               \frac{dF}{ds} = -\frac{\partial \widehat{H}}{\partial z}
                              =  \frac{1}{2}p. \nonumber
          \ee
       The aim of this note is to prove the following theorem using the differential Galois
       approach.

       \begin{th1}
          The Hamiltonian system (\ref{eqh}), corresponding to
          ${\rm P}_{{\rm II}}$, does not have second meromorphic first integral for
          $\alpha = n \in \ZZ.$
      \end{th1}

        \noindent
         {\bf Remark.} We take these values of the parameter
        $\alpha = n \in \ZZ$ because it is known that for every
        $\alpha \in \ZZ,\, {\rm P}_{{\rm II}}$ has unique rational solution. We
        need simple single-valued solutions to study variational
        equations along them.

        In section 2 we recall some notions and fact about
        Differential Galois Theory and Morales - Ramis theorem.
        The proof of Theorem 1 will be given in section 3.

        Before ending this section we  summarize some facts about
        ${\rm P}_{{\rm II}}$ which are useful for our purpose.

         Suppose $\alpha = n \in \ZZ$
        in (\ref{pe}) then
          \be
            \label{ope}
            w(z, -n) = -w(z, n)
          \ee
        which reduces our considerations
        in Theorem 1 to the cases $n \in\ \NN_0 := 0, 1, 2, \ldots$
        (see section 3).
        Next, the B\"acklund transformation
            \be
              w(z, n + 1) = -w(z, n)-
                      \frac{2n + 1}{2w^2(z, n) + 2w'(z, n) + z}\,,
            \ee
        generates the hierarchy of rational solutions of ${\rm P}_{{\rm II}}$
        from the trivial solution $w(z, 0) = 0.$ We list only few of them
        $w(z, 1) = -\frac{1}{z}, w(z, 2) = \frac{1}{z} - \frac{3 z^2}{z^3 + 4},
        w(z, 3) = \frac{3 z^2}{z^3 + 4} - \frac{6 z^5 + 60 z^2}{z^6 + 20 z^3 - 80}, \ldots$.

       Following $[^5]$ we define the Vorobev-Yablonski polynomials $Q_n(z)$ by
        the recursion relation
            \be
               Q_{n + 1}\,Q_{n - 1} = zQ^2_n + 4 (Q'_n)^2 - 4 Q_n\,Q''_n
            \ee
       with $Q_0(z) = Q_1(z) = 1.$
       It turns out that $Q_n(z)$ are monic polynomials of degree
       $\frac{1}{2}n(n - 1).$ Then the rational function
           \be
             \label{pol}
               w(z, n) = \frac{d}{dz}\Big\{\ln\Big[\frac{Q_n(z)}{Q_{n + 1}(z)}\Big]\Big\}
                   = \frac{Q'_n(z)}{Q_n(z)} - \frac{Q'_{n + 1}(z)}{Q_{n + 1}(z)}
           \ee
       satisfies ${\rm P}_{{\rm II}}$.
       We need also the following two theorems.
      \begin{th1} $[^6]$
           For every positive integer $n$, the polynomial $Q_n(z)$ has simple
           roots.
      \end{th1}
     \begin{th1} $[^6]$
            For every positive integer $n$, the polynomials $Q_n(z)$ and
           $Q_{n + 1}(z)$ do not have  common roots.
\end{th1}

\section{Theory }

    The aim of this section is to formulate basic definitions and facts from
    Morales - Ramis theory following $[^{3,2}]$.
    Let us consider a $2n$--\,dimensional complex analytic
    manifold $M$ and a holomorphic Hamiltonian system $X_H$ on $M$
        \be
           \label{tham}
            \frac{d}{dt} z(t) = X_H (z).
        \ee
    Let $z = z(t)$ be a  non--\,equilibrium solution of (\ref{tham}).
    The phase curve $\Gamma$
    is the connected Riemann surface corresponding to $z = z(t).$
    The variational  equations $(VE)$ along $z = z(t)$ have the form
         \be
          \label{tve}
           \frac{d}{dt}\, \xi = X'_H(z(t))\,\xi,\qquad\
              \xi\in T_{\Gamma}M.
         \ee
    Following Ziglin $[^7]$ we can always reduce the order of this system by one. Let
    $N:=T_{\Gamma} M/T\Gamma$ be the normal bundle of $\Gamma$ and
    $\pi:T_{\Gamma} M\rightarrow N$ be the projection. Then the
    system (\ref{tve}) induces the following system on $N$
          \be
           \label{tnve}
             \frac{d}{dt}\,\eta = \pi_{\star}(X'_H(z(t))\pi^{-1}\eta),
             \qquad \eta \in N
          \ee
   which is called the normal variational equation $(NVE)$
   along $\Gamma.$ We shall complete the Riemann surface $\Gamma$
   with some equilibrium points and (possibly) points at infinity,
   in such a way, that the coefficients of the (VE) and of the
   (NVE) are meromorphic on this extended Riemann surface
   $\overline{\Gamma}\supset \Gamma$  (see $[^{2,3}]$).

   Morales-Ruiz and Ramis formulate a criterion for
   non-integrability for Hamiltonian systems in terms of the
   properties of the differential Galois group $G$ of the normal variational
   equations $[^3]$.
\begin{th1} $[^2]$
           Assume there are $n$  first integrals of $X_H$
           which are meromorphic in a neighborhood
           of  $\overline{\Gamma}$, in particular meromorphic at infinity and independent
           in a neighborhood
           of  $\overline{\Gamma}$ (not necessarily on $\overline{\Gamma}$ itself). Then the identity component
           of the Galois group of the (VE) (resp.  (NVE)) is
           abelian.
\end{th1}

   In the applications first a non-equilibrium particular solution is
   selected. Next, we calculate
   the $VE$ and $NVE$. And finally, we have to check if
   the identity component of differential Galois group $ G^0$
   of obtained $NVE$ is Abelian. If it is not, then
   the Hamiltonian system $X_H$ is not integrable.

    Computing the differential Galois group is usually difficult but in
    the case of second order equation with rational coefficients
    there exists an  efficient algorithm - the Kovacic algorithm. We briefly recall it -
     see $[^{9,10}]$ for more detailed description.

   Any second order differential equation with rational
   coefficients can be reduced to the form
     \be
       \label{requ}
          y'' = r y.
     \ee
   The logarithmic derivative $\omega := \frac{y'}{y}$ of solution
   $y$ of (\ref{requ}) satisfies the Riccati equation
     \be
       \label{ric}
         \omega' + \omega^2 = r
     \ee
  and according to Lie - Kolchin's theorem, equation (\ref{requ})
  has a Liouvilian solution (that can be expressed via exponentials,
  integrals and algebraic functions) if and only if the corresponding
  Riccati equation (\ref{ric}) has an algebraic solution
  (for definitions, details and proof related to differential
  algebra see $[^{8,3,9}]$).
  Moreover, the degree $m$ of the minimal polynomial for this
  algebraic solution is one of the following numbers
    $ L_{max}:= \big\{1, 2, 4, 6, 12\big\}.$
  The differential Galois group $G$ of (\ref{requ}) is
  an algebraic subgroup of ${\rm SL}(2, \CC)$ and
   has one of the following forms

      \begin{enumerate}

            \vspace{-1ex}
            \item
     Case 1: $ G$ is triangularisable; in this case
     equation (\ref{requ}) is reducible and has a solution of the
     form
     $y = \exp \int \omega$, where $\omega \in \CC(z)$, i.e., Riccati
     equation (\ref{ric}) has a rational solution $(m = 1)$.

            \vspace{-1ex}
            \item
     Case 2: $ G$ is ''imprimitive'' i.e. conjugate to a subgroup of
     $D \bigcup          \left(
           \begin{array}{cc}
              0 & 1  \\[0.6ex]
               -1      & 0
           \end{array}
                   \right) D$, where $D$ is the diagonal subgroup in ${\rm SL}(2, \CC)$ ;
      in this case
     equation (\ref{requ}) has a solution of the
     form
     $y = \exp \int \omega$, where $\omega$ is algebraic over $\CC(z)$
     of degree 2, i.e., Riccati
     equation (\ref{ric}) has an algebraic solution of degree $m = 2$.

            \vspace{-1ex}
            \item
     Case 3: $ G$ is finite and cases 1 and 2 do not hold; for this case
     all solutions of equation (\ref{requ}) are algebraic and
     Riccati equation (\ref{ric}) has an algebraic solution of degree
     $m \in \big\{4, 6, 12\big\}$.

            \vspace{-1ex}
            \item
     Case 4: $ G = {\rm SL}(2, \CC)$ and
     equation (\ref{requ}) has no Liouvillian solution, i.e., Riccati
     equation (\ref{ric}) has no algebraic solution.

     \end{enumerate}

\section{Proof of Theorem 1}

     Consider first the case $\alpha = n \in \NN$. The Hamiltonian
     system (\ref{eqh}) possesses a particular solution of the
     kind
       \be
         \label{psol}
            q(s) = w(s, n), \qquad
                 p(s) = w'(s, n) + w^2(s, n) + \frac{1}{2}\,z(s) , \\
            z(s) = s , \qquad
                 F(s) = \frac{1}{2}\int p\,ds \nonumber
       \ee
     where $w(s, n)$ are already defined rational functions (\ref{pol}).
     Note that $q(s), p(s), z(s), F(s)$ are rational functions,
     i.e., the particular solution is single valued.
     The phase curve $\Gamma$ of the solution (\ref{psol})
     rationally parameterizes the $\widehat{H} = 0$ and naturally
     extends to $\overline{\Gamma}:= \CP$.

     The variational equations $(VE)$ along $\overline{\Gamma}$
     are
      \ben
       \frac{d}{ds} \left(
                          \begin{array}{c}
                                \xi_1 \\[0.6ex]
                                \eta_1 \\[0.6ex]
                                \xi_2 \\[0.6ex]
                                \eta_2
                          \end{array}
                    \right)=
                    \left(
           \begin{array}{ccrc}
              -2w(s, n)   & 1                 &-\frac{1}{2}   &0
              \\[0.6ex]
               2p(s)      & 2w(s, n)          &0              &0  \\
               [0.6ex]
               0          & 0                 &0              &0  \\
               [0.6ex]
               0          &\frac{1}{2}        &0              &0
           \end{array}
                   \right)
                   \left(
                          \begin{array}{c}
                                \xi_1 \\[0.6ex]
                                \eta_1 \\[0.6ex]
                                \xi_2 \\[0.6ex]
                                \eta_2
                          \end{array}
                    \right)
     \een
    and we can take the upper left block as normal variational
    equations $(NVE)$ (as $z = s$ we substitute $\frac{d}{ds}$
    with $\frac{d}{dz}$)
        \be
          \label{nve}
            \frac{d}{dz}\,\xi_1  &=& -2 w(z, n)\,\xi_1 + \eta_1
            \\[0.3ex]
            \frac{d}{dz}\,\eta_1 &=& 2 p(z)\,\xi_1 + 2 w(z, n)\,\eta_1
            \nonumber
        \ee
    In order to study the differential Galois group of (\ref{nve})
    we reduce it to a single second order equation
       \be
          \label{rnve}
           \xi'' = (6\,w^2(z, n) + z)\, \xi.
       \ee
    Now we apply the Kovacic algorithm to the equation (\ref{rnve}).
    Denote
              $$r(z):=
              6w^2(z, n)+z=6\Big[\frac{Q^{'2}_n(z)}{Q^2_n(z)}+
                             \frac{Q^{'2}_{n + 1}(z)}{Q^2_{n+1}(z)}-
                            2\frac{Q'_n(z)Q'_{n + 1}(z)}{Q_n(z)Q_{n + 1}(z)}\Big] + z
                            := \frac{R(z)}{S(z)}$$
   where
      $
         \deg\,S(z) = \deg\,Q^2_n(z)\,Q^2_{n+1}(z) = 2n^2,
         \deg\,R(z) = \deg\,z\,Q^2_n(z)\,Q^2_{n+1}(z) = 2n^2 + 1.
      $

   Let $Y'$ be the set of zeroes of $S(z)$
      $$Y' = \big\{c\in\CC\,|\,S(c) = 0\big\}
          = \big\{c \in \CC\,|\,Q_n(c) = 0,\,Q_{n + 1}(c) = 0\big\}$$
   and $Y = Y'\cup \big\{\infty\big\}$.
   The order of $c$ denoted by $o(c)$ is the multiplicity of $c$ as a root of
   $S(z)$ and due to Theorem 2 and Theorem 3, $o(c) = 2$ for every
   $c\in Y'$. The order of infinity is
     $o(\infty) = \max(0, 4 + \deg\,R(z) - \deg\,S(z)) = 5.$

  Define
    $m^+ := \max_{c\in Y} o(c)= o(\infty) = 5$
  and $Y_i = \big\{c\in Y\,|\,o(c) = i \big\}$. In our case
    $$Y_2 := \big\{c\in Y\,|\,o(c) = 2 \big\}
          = \big\{c\in\CC\,|\,S(c) = 0 \big\},\qquad
          Y_5 := \big\{\infty\big\}.$$
  Then
    $\gamma_5 := card\,Y_5 = 1,\quad
      \gamma_2:= card\,Y_2 = \frac{1}{2}n(n-1) + \frac{1}{2}n(n+1) = n^2,$

  so we compute
     \ben        \gamma\,    &:=&  \,\gamma_2 +
     \sum_{\substack{{\rm odd \, k}  \\
                       3\leq k\leq m^+}}\gamma_k = n^2 + \gamma_5 = n^2 + 1 \geq 2.
     \een
  Define   $L' \subset L_{max}:L' = \big\{2\big\}$ and since
  $m^+ = 5 > 2$ then $L \equiv L'= \big\{2\big\}$ or $m = 2,\,h(2) = 2$,
  i.e., for the Galois group of (\ref{rnve}) we have only cases 2
  or 4.

 For every $c\in Y_2$ we need the Laurent series expansion of
 $r$ around $c$
    $$r = \frac{\alpha_c}{(z - c)^2} + \frac{\beta_c}{z-c}+O(1).$$
 Recall that $c = z_j,\quad j = 1,\ldots, n^2$ are simple roots of
 $Q_n(z)$ and $Q_{n + 1}(z)$. Obviously $\alpha_{z_j} = 6$. It remains
 to compute
  $\Delta_{z_j} = \sqrt {1 + 4 \alpha_{z_j}} = 5.$

  We proceed with the second step of the Kovacic algorithm
  defining $(m = 2)$
   \ben
     E_{z_j}    &=&
     \big\{2-(2-2j)\Delta_{z_j}\,|\,j = 0, 1, 2\big\}\cap \ZZ
            = \big\{-8, 2, 12\big\}\qquad
            \textrm{and} \\
     E_{\infty} &=&
     \big\{o(\infty)\big\} = 5.
   \een
   Finally, in the third step, we note that there is no such
   $e = \big\{e_c\big\}_{c\in Y}$ in the Cartesian product
   $E := \prod_{c\in Y}E_c$ for which
     $d(e) = 2 - \frac{1}{2}\sum_{c\in Y} e_c \in \NN_0$
   since there is only one odd number in the sum - $e_{\infty}$ and
   the others are even. Then the Galois group of (\ref{rnve}) is
   ${\rm SL}(2, \CC)$, i.e., not Abelian.

   It remains to deal with the case $\alpha = 0.$ The particular
   solution of (\ref{eqh}) in this case is
     \be
       \label{npsol}
         q = 0, \quad p = \frac{1}{2}\,s,\quad z = s,\quad
         F = \frac{1}{2}\int p\,ds.
     \ee
  In the same manner, the normal variational equation is obtained
  to be $         \xi'' = z\,\xi.      $
  This is  Airy equation and its Galois group is
  ${\rm SL}(2, \CC)$ $[^9]$.
  Note that  from (\ref{ope}) and that $r(z)$ is an even function with respect to $n$,
    it  follows that above consideration  is
  valid for $\alpha = n \in \ZZ.$
  Hence, from Morales-Ramis Theorem (Theorem 4.) the Hamiltonian
  system is not completely integrable with meromorphic first
  integrals.  This ends the proof.

  {\bf Remark.} Actually, for these $\alpha$ we prove that the
  system (\ref{eqh}) does not possess another rational first
  integral except $\widehat{H}$--the meromorphic functions on
  $\overline{\Gamma} = \CP$ are exactly the rational ones.

  As seen from the proof, equation (\ref{rnve}) and
   Airy equation are not Fuchsian (the point $z = \infty$ is an irregular
  singular point). Such points are considered the main source of
  non - integrability because of existence of a non-trivial Stokes
  phenomenon at them (see $[^3]$ for example ).

  Finally, we recall that in $[^{11}]$  it is shown that ${\rm P}_{{\rm II}}$
(with slightly different Hamiltonian) has nonrational
(possibly multi - valued) second integral.

{\bf Acknowledgements}. This work is partially supported by grant
MM 1504/05 with NSF Bulgaria.

\medskip

 REFERENCES

 $[^1]$ Ablowitz M., P. Clarkcon, Solitons, Nonlinear Evolution Equations and
                Inverse Scattering. Cambridge, Cambridge University Press, 1991.
  $[^2]$ Morales-Ruiz J. Regul. Chaotic Dyn., {\bf 5}, 2000, 251.
  $[^3]$   Morales-Ruiz J. Differential Galois Theory and
                Non-Integrability of Hamiltonian Systems. Progress in
                Mathematics {\bf vol 179}, Basel: Birkh\"au\-ser, 1999.
 $[^4]$    Arnold V. Mathematical methods in classical mechanics. Berlin,
                Springer-Verlag, 1978.
 $[^5]$    Vorobev A. Differ. Uravn., {\bf 1} 1965, 58.
 $[^6]$  Fukutani S., K. Okamoto, H. Umemura. Nagoya Math.
                J., {\bf 159}, 2000, 179.
$[^7]$    Ziglin S. I, Funct. Anal Appl. {\bf 16}, 1982, 181;
                 II, Funct. Anal Appl. {\bf 17}, 1983, 6.
$[^8]$   Put M. van der, M. Singer.
                 Galois Theory of linear differential equations.
                Grundehren der Mathematischen Wissenschaften, {\bf 328}, Berlin,
                Springer, 2003.
 $[^9]$    Kovacic J. J. Symb. Comp., 1986, 3.
$[^{10}]$    Maciejewski A., J. Strelcyn, M. Szydlowski. J. Math. Phys.,
                 {\bf 42}, 2001, 1728.
$[^{11}]$   Flaschka H., A. Newell. Comm. Math. Phys., {\bf 76}, 1980, 65.

\begin{flushright}
Department of Mathematics and Informatics,\\
 Sofia University \\
5 James Bouchier blvd., \\
1164 Sofia, Bulgaria \\
e-mail: cveti@fmi.uni-sofia.bg
\end{flushright}

\end{document}